# The Second Postulate of Euclid and the Hyperbolic Geometry

Yuriy N. Zayko
Department of Applied Informatics, Faculty of Public Administration,
Russian Presidential Academy of National Economy and
Public Administration, Stolypin Volga Region Institute, Russia, Saratov.
zyrnick@rambler.ru

**Abstracts**. The article deals with the connection between the second postulate of Euclid and non-Euclidean geometry. It is shown that the violation of the second postulate of Euclid inevitably leads to hyperbolic geometry. This eliminates misunderstandings about the sums of some divergent series. The connection between hyperbolic geometry and relativistic computations is noted.

1. Introduction

Historically, the appearance of non-Euclidean geometry is associated with the realization of the possibility of not fulfilling the fifth principle (postulate) of Euclid about parallel lines. Geometry introduced in the writings of Lobachevsky and Boyayi instead of the fifth postulate of Euclid takes the opposite one and is just as consistent as the Euclidean geometry. It was named hyperbolic [1]. The selection of the fifth postulate laid the foundation for the accepted division of geometries into absolute one based on the first four postulates of Euclid, Euclidean geometry, in which, in addition to the first four, the fifth postulate is added and the hyperbolic geometry already mentioned.

The distinct feature of the fifth postulate from the others was stressed long before the appearance of non-Euclidean geometry. The rest of the postulates did not cause such increased attention, and, especially, doubts about their fairness, which seems rather strange, since there are numerous examples of violation of at least one of the rest postulates of Euclid - the second postulate.

As shown in this paper, this indicates a possible deviation from the Euclidean geometry in the rejection of the second postulate, which should also be taken into account in mathematical and other studies. To bind it with the violation of the fifth postulate is not possible, therefore, it is necessary to carefully study the second postulate, what is done below.

2. The second postulate of Euclid

Let us give below one of the formulations of the second postulate [1]

*A finite straight line may be extended continuously in a straight line*

Like any statement expressed in verbal form, it differs in ambiguity and admits numerous variants. For example, in [2] the word "continuously" is replaced by "unlimited".

In order to eliminate this inaccuracy, we resort to a technique common in the mathematics - the modeling of the statements of one region by the means of another [1]. In this case, the second postulate in the language of arithmetic is equivalent to the following

*A sum of an infinite divergent series of positive numbers is equal to infinity*

In accordance with it, for example, there should be

$$1+1+1+1+\ldots+1+\ldots = \infty$$
$$1+2+3+4+\ldots+n+\ldots = \infty \tag{1}$$

However, this contradicts known facts [3]. The sum of the first series is equal $\zeta(0) = -0.5$, and of the second, $\zeta(-1) = -1/12 = -0.833...$, where $\zeta(s)$-is the Riemann zeta-function [3]

$$\zeta(s) = \sum_{n=1}^{\infty} n^{-s}, s = u + iv \tag{2}$$

$n$ – is an integer, $u$ and $v$ –are real numbers. There are a lot of such examples.

The very fact of the finiteness of the sum of a divergent series does not raise questions since its meaning is different than for the sum of convergent series [4]. Divergent series are used in various fields of science, primarily in physics. In particular, the sum of the second series in (1) underlies many results of string theory [5]. However, until recently no one paid any attention to the fact that this and similar results violate the second Euclidean postulate and the consequences of this fact.

3. Calculation of the sum of a divergent series

Recall that the sums of divergent series are not computed, since direct computation is usually impossible, but is determined either indirectly, for example, as in the case of the series for $\zeta(-1)$ - by analytic continuation of the zeta-function, or by using various summation methods [4]. The first attempt to calculate the sum of this series was undertaken in the author's work [6]. It turned out that for this it is necessary to introduce a metric on the numerical axis

$$ds^2 = \left(1+\frac{x}{x_c}\right)c^2 dt^2 - \left(1+\frac{x}{x_c}\right)^{-1} dx^2 \qquad (3)$$

where $s$ – is an interval; $x, t$ – are coordinate on the numerical axis and time; $x_c$ – some characteristic value for the given metric [6]. The calculation is realized as the motion of a material particle along the numerical axis according to with the relativistic equations of motion written for the metric (3). Factually, it is talking about the relativistic Turing machine (a relativistic supercomputer in the terminology of [7]) whose tape is represented by the numerical axis with the metric (3), and the role of moving head is played by the above-mentioned particle. Such a supercomputer is able to solve problems that are not computable in the traditional sense, in this case, to calculate the sum of a divergent series. The accuracy attained in [6] corresponds to an error in 3,5%.

4. Nature of Geometry

Now it is understood the perplexity that caused by the received result - the magnitude of the sum of $\zeta(-1)$ for a number of representatives of the academic community [8]. It was caused by an attempt to interpret it in accordance with the so-called common sense, i.e. in fact, an attempt to comprehend it within the framework of Euclidean geometry. In any case, no deviations from it, not related to the violation of the fifth postulate, no one expected. Nevertheless, this is so.

The character of the resulting geometry can be easily determined if one tries to solve the problem of embedding a one-dimensional manifold with the metric (3) into a manifold of higher dimension - the plane. As will be shown below, in whole, for $-\infty < x < \infty$ this is possible, only for a plane with hyperbolic geometry. Before tackling directly the solution of this problem, we recall that the awareness of the impossibility of a smooth embedding of a subspace into a Euclidean space of higher dimension in its time contributed to the discovery of non-Euclidean geometry [9].

Thus, consider a two-dimensional plane with a spatial metric defined on it

$$dl^2 = \pm dx^2 \pm dy^2 = \left[\pm 1 \pm \left(\frac{dy}{dx}\right)^2\right]dx^2 \qquad (4)$$

Comparing with (3), find that the embedding leads to the equation for *y(x)*

$$\pm 1 \pm \left(\frac{dy}{dx}\right)^2 = \frac{1}{1 + \frac{x}{x_c}} \qquad (5)$$

which is easily integrated. Below are the results, according to which there are three areas, each of which is characterized by its own function *y(x)*

I. $\quad -\infty < x < -x_c; dl^2 = dx^2 - dy^2;$
$\quad\quad y(x) = \pm x_c \left(arcch\sqrt{-z} + \sqrt{(1+z)z}\right)$

II. $\quad -x_c < x < +\infty; dl^2 = dy^2 - dx^2;$
$\quad\quad y(x) = \pm x_c \left(arcsh\sqrt{1+z}\right) + \sqrt{(1+z)(2+z)}$ $\qquad (6)$

III. $-x_c < x < 0; dl^2 = dx^2 + dy^2;$
$y(x) = \pm x_c \left(\arcsin\sqrt{-z} - \sqrt{-z(1+z)}\right)$

where $z = x/x_c$. This is illustrated in Figure 1, where some branches of the embedding (6) are shown

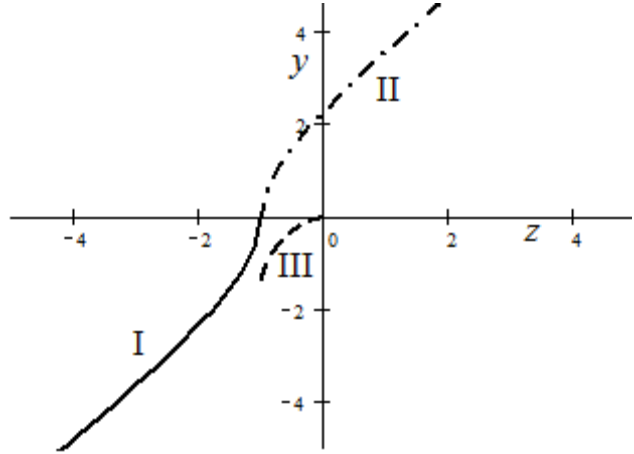

Fig. 1. Areas of embedding of a one-dimensional manifold with metric (3) into a two-dimensional plane. The Roman numerals next to each branch of $y(x)$ refer to each corresponding region: I-solid line, II-dash-dotted, III- dashed.

From these results, it follows that one can completely embed a one-dimensional manifold with the metric (3) only in a two-dimensional plane with hyperbolic geometry (branches I and II). Into a plane domain with Euclidean geometry, only part of it corresponding to the domain III can be embedded.

5. Discussions.

The fact that the violation of the second postulate of Euclid leads to non-Euclidean geometry, restores the "symmetry" between postulates, abolishing the "monopoly" of the fifth postulate, which lasted about two hundred years. This is a comparatively short period of time if one compares it with the time during which Euclidean geometry dominated. It should be noted that the rejection of both postulates leads to the same consequences.

It is now difficult to foresee the other consequences to which this result will lead. Right now it helps to eliminate misunderstandings about the sums of some divergent series mentioned above. As for the possibilities in the field of computational mathematics, and in particular the creation of relativistic supercomputers, it is worth noting the undoubted advantages of this approach in comparison with the traditional one where it is proposed to use Kerr-Newman black holes for the realization of relativistic calculations [7].

Calculation of the $\zeta(-1)$, made in [6], as already noted above, has insufficient accuracy. If one considers it as another confirmation of the general theory of relativity, then it is inferior in accuracy in comparison with all the other (experimental). However, in view of their fewness, it should not be discounted.

It must be said that the calculation method of $\zeta(-1)$ used in [6] is physical, and based on the fact that the formula for the partial sums of the series (1) coincides with the expression for the distance traversed by a particle moving with a constant (nonrelativistic) acceleration. Strictly speaking, it was this analogy that made it possible to obtain expressions for the metric (3). It is not yet clear how this result can be obtained from purely geometric considerations. One can only hope that the establishment of a closer relationship between the two approaches will allow for greater progress in relativistic calculations and, in particular, to improve their accuracy.

6. Conclusion

The paper shows that geometry in which the second postulate of Euclid is not satisfied is hyperbolic. This allows eliminating the misunderstandings associated with the calculation of the sums of some divergent series. A metric is given on a numerical axis in which the specified calculations can be performed. The embedding of a numerical axis with a given metric into a two-dimensional plane was performed and it is shown that the condition for smooth embedding of the entire numerical axis is the hyperbolic geometry on the plane. Various applications of the obtained results are considered.

References.

1. H.S.M. Coxeter, F.R.S. Introduction to Geometry, New York, London, John Wiley&Sons, Inc.
2. S.N. Bychkov, E.A. Zaitsev, Mathematics in world culture. Tutorial, Moscow: Russian State Humanitarian Univ. Publishers, 2006 (Russian).
3. E. Janke, F. Emde, F. Lösch, Tafeln Höherer Funktionen, B.G. Teubner Verlagsgesellschaft, Stuttgart, 1960.
4. G.H. Hardy, Divergent series, Oxford, 1949.
5. B. Zwiebach, A First Course in String Theory, 2-nd Ed., MIT, 2009.
6. Y.N. Zayko, The Geometric Interpretation of Some Mathematical Expressions Containing the Riemann $\zeta$-Function, Mathematics Letters, 2016; 2(6): 42-46.
7. H. Andréka, I. Németi, P. Németi, General relativistic hypercomputing and foundation of mathematics, Natural Computing, 2009, V. 8, № 3, pp. 499–516.
8. D. Berman, M. Freiberger, Infinity or -1/12?, + plus magazine, Feb. 18, 2014, http://plus.maths.org/content/infinity-or-just-112
9. S. Weinberg, Gravitation and Cosmology. Principles and Applications of the General Theory of Relativity, MTI, John Wiley&Sons, 1972.